# Determining Optimal Lot Size, Reorder Point, and Quality Features for a Food Item in a Cold Warehouse: Data-Driven Optimization Approach

Atena Karimi, Kharazmi University, Department of Industrial Engineering, Tehran, Iran

atenakarimi2000@yahoo.com

Omid Ghorbani, Kharazmi University, Department of Industrial Engineering, Tehran, Iran

omidghorbaniai2001@gmail.com

Reza Tashakkori, Kharazmi University, Department of Industrial Engineering, Tehran, Iran

rezatashakkory@gmail.com

*Seyed Hamid Reza Pasandideh, Kharazmi University, Department of Industrial Engineering, Faculty of Engineering, Kharazmi University, Tehran, Iran (Corresponding Author)

OrcID: 0000-0003-2468-8682

shr_pasandideh@khu.ac.ir

**Abstract**

We propose a nonlinear optimization model for determining the optimum lot size and reorder point for a food item distributed through a cold warehouse as well as the optimum quality features, namely temperature, humidity, packaging type, and level of environmental conditions. The item's quality is estimated based on the features mentioned earlier, and then it is used as a constraint in the optimization process. An assumption was made that the inventory is managed under a continuous review policy and the warehouse has limited space. The model seeks to minimize the annual total cost of managing the warehouse. The model will be a nonlinear mixed programming one, which is solved by Pyomo as a leading library in Python language programming. Numerical examples are used to demonstrate the use of the model and, through sensitivity analysis, develop insights into the operation of cold warehouses. This sensitive analysis opens the doors to managerial insight from which managers and policymakers can highly benefit.

## 1 Introduction

The changing lifestyle of consumers and their increased awareness of the health benefits of consuming fresh food have increased their demand. This growth in demand has elevated the importance of cost-effective management of the fresh food supply chain. Since the variety of fresh foods is extensive and perishable, managing their supply chain is costly and challenging. Fresh food items have a natural harvest-to-table duration beyond which they are deemed perished and unsuitable for human consumption. The natural harvest-to-table duration of an item is the elapsed time between its harvest and its consumption by humans without any intervention to prevent the growth of adversarial microbial colonies affecting its safety. The natural harvest-to-table duration of a food item depends on its chemical and protein composition and its respiration rate. Since these characteristics vary widely among food items, their natural harvest-to-table durations are diverse. The shortness of the natural harvest-to-table duration is an impediment to satisfying the consumer demand for the availability of high-quality fresh food. Therefore, to mitigate this impediment, suppliers, with the help of food scientists, have relied on cold storage to extend the natural harvest-to-table duration of fresh food. Cold temperatures slow or impede the growth of microbial communities; and in general, colder temperatures result in slower growth of microbial communities. In general, managing a cold supply chain is costlier and more challenging than managing conventional supply chains. The increase in cost is primarily due to the energy cost of refrigeration, special food handling requirements, and food safety regulations. Additional challenges stem from additional constraints imposed on inventory management and logistics. This paper focuses on the warehouse inventory management of fresh food items, the first stage in the cold supply chain.

## 2 Literature Review

In this paper, the studies which are relevant to the research are classified into three parts, each of which tries to shed light on different aspects.

## 2.1 Cold Supply Chains

In cold chain logistics, the importance of efficient management for ensuring the integrity and safety of perishable goods during transportation is widely recognized. Accordingly, [1] have developed a compact supervisory system specifically designed for cold chain logistics organizers to track and identify improper practices for refrigerated transportation based solely on temperature data. It is important to note that cold supply chains contribute significantly to greenhouse gas emissions, and this issue urgently needs to be addressed. Thus, [2] proposed a mixed-integer concave minimization problem for cold supply chain design. As a result of this approach, total costs are minimized while global warming is reduced. Similarly, In a study conducted by Goodarzian et.al a bi-objective mixed-integer non-linear programming problem model was created to optimize the citrus fruit supply chain by simultaneously minimizing costs and CO2 emissions. To solve this complex optimization problem, researchers used two techniques encompassing the Augmented Epsilon-constraint method and Multi-objective optimization algorithms [3]. Regarding implementing "green" practices that can enhance the sustainability of supply chains, Ali et al. applied a three-phase approach, merging contingency theory and the triple bottom line concept, to assess sustainability in supply chains. Focusing on frozen food warehouses in Saudi Arabia, they used fuzzy Delphi and the best-worst method to identify and prioritize sustainable practices [4]. Also as a result of the depreciation of perishable goods over time, the distribution of products efficiently plays a crucial role in the successful operation of supply chains. Hence, [5] introduced a novel mixed-integer mathematical formulation for optimizing truck scheduling at a cold-chain cross-docking terminal to minimize the total cost incurred during truck service. Moreover, a customized Evolutionary Algorithm has been developed to tackle the complexity of the model. It is worth mentioning that [6] addressed the overlooked government policy considerations in reducing energy usage and emissions within cold chain logistics enterprises. To tackle this issue, they proposed a decision-making model employing a method called chaotic particle swarm optimization (CPSO). Moreover, study of Liu et al. [7] analyzed the development of cold chain logistics in China, focusing on insufficient terminal logistics distribution capacity in regional centers. As part of this study, a joint distribution node was proposed to improve terminal logistics distribution capacity by transporting goods using semi-trailers and integrating cold chain logistics resources effectively. Lastly, A study conducted by Win Wang and Lindu Zhao revealed that collaborative cold chain investments and pricing were superior strategies for maximizing profits in the fresh food supply chain. Cold chain investments are positively influenced by consumer food safety awareness in safety-sensitive markets [8].

## 2.2 Applications of Data Science in Supply Chains

As Machine Learning (ML) state-of-the-art approaches and methods have recently gained popularity as advantageous tools in Predictive Maintenance (PdM) applications, [9] found it helpful and applicable to use proper ML tools for designing PdM applications to ensure their success. It is essential to point out that as advanced monitoring and diagnosis technologies have become more widely available, ML has become increasingly important. Relying on the fact that most companies find it difficult and costly to implement ML-based and PdM based on ML is the most appropriate strategy, the study by [10] defined a cost-oriented analysis and provided a mathematical model. This model gets assistance from different ML tools and classifiers to evaluate its accuracy and performance. To delve deeper into the advantages of ML prediction methods, which outperform traditional approaches in terms of accuracy and their ability to handle large and complex datasets, a study by [11] has been conducted. This implies that they implemented a new two-stage solution approach for supplier selection and order allocation planning, integrating a prediction procedure with an optimization model. As regards, managing supply chain risks in recent years has become increasingly significant, [12] introduced a framework for predicting supply chain risks, which combines data-driven Artificial Intelligence (AI) techniques and the expertise of supply chain professionals. Following this, the framework is applied to the problem of predicting delivery delays in a multi-tier manufacturing supply chain. The objective was to assess the trade-off between prediction performance and interpretability by employing various ML methods. Ultimately, in another momentous study conducted by [13], a new prediction model was developed to accurately forecast time-varying volatile processes without relying on linearity assumptions in time series data. As a matter of fact, the researchers formulated the prediction problem as an optimization problem, employing a data-driven approach. The effectiveness of this approach was evaluated by utilizing multiple datasets obtained from real-life applications. The numerical results derived from their model convincingly demonstrated its competitive prediction capabilities, providing compelling evidence of its effectiveness. Last but not least, [14] developed a stochastic optimization model for apple supply chain management, integrating supplier selection and cold storage inventory decisions. Through a real-world case study, they demonstrate significant cost reductions from quality-focused optimization under uncertainty. The study provides an effective model to support agribusiness production planning and uncertainty mitigation over harvesting seasons. The study of Tufano et al. explores the potential of underutilized data from warehouse management systems, investigating its broader costly implications beyond its original purpose by employing some ML classifiers to predict crucial aspects of storage systems and also assessing how more data availability affects higher prediction accuracy to emphasize the managerial relevance of this data-driven approach [15]. The results of Alsina et al. research showcased the superior effectiveness of the data-driven approach combined with various ML classifiers, particularly random forest, highlighting its consistently high accuracy in predicting component reliability and its substantial improvement over traditional methods, especially with larger time-to-failure datasets [16].

### 2.3 Applications of Data Science in Optimization

The field of supply chain management is increasingly adopting ML techniques to address complex challenges and optimize decision-making processes. This includes the realm of cold chain management whereby [17] identified an unexplored area of detecting cold chain breaks and proposed using anomaly detection methods in time series data. By harnessing the power of ML, operators can deepen their understanding of cold chain disruptions and receive timely alerts to prevent further damage. In response to the escalating traffic congestion in megacities, [18] conducted a data-driven study using an enhanced weighted K-means clustering algorithm on actual traffic congestion index data from Beijing. They aimed to identify and analyze traffic congestion patterns to inform tailored optimization and control strategies, ultimately mitigating congestion and enhancing traffic equilibrium with the help of ML tools. Data mining techniques impact the progress of retail marketing, especially with the ever-growing amount of customer data. The study of [19], presented an ML algorithm on a retail marketing dataset that integrates Multi Variant K means clustering and pattern mining to forecast customer interests, explicitly focusing on assisting E-Commerce systems. In [20] study, researchers analyzed the factors that impact consumer's willingness to pay for cold chain logistics using survey data from 711 consumers across nine cities in China. [21] employed a data-driven optimization approach to design a reverse logistics network for a disposable product recycling management system. Similarly, the study of [22] introduces two innovative data-driven game-theoretic methods leveraging ML to enhance inventory management and financial choices for sellers receiving online platform financial assistance. The results demonstrate the advantages of integrating ML with quantile regression for sellers and lenders on 199,390 weekly sales records of retail transaction data. Notably, in demand forecasting for supply chain performance enhancement, [23] compared traditional forecasting methods with ML-based techniques. Their analysis revealed significant performance differences and highlighted the potential of hybrid forecasting methods that incorporate ML models like ARIMAX and Neural Networks. Consequently, [24] explored the potential of sensor data in managing the chilled food chain, using a prototype tracking tool to predict the remaining shelf-life. Their study showcased the benefits of sensor data for pricing decisions and performance improvement, underscoring the impact of sensor data-driven supply chain management. To further elucidate, ML applications have played a crucial role in bolstering sustainability within agricultural supply chains. In line with this, [25] conducted a comprehensive systematic review that showcased the potential advantages of employing ML algorithms across various stages of agricultural supply chains (ASCs). Furthermore, [26] study emphasizes the increasing importance of improving microgrid system utilization and reducing dependence on traditional energy sources in warehouse operations, addressing operational challenges through a mixed linear programming model and a data-driven variable neighborhood search algorithm for task optimization. Last but not least, in the context of vaccine supply chains, (Lin, Zhao, and Lev 2020) examined the decision-making process surrounding transportation choices and retailer roles. Their study emphasized the influence of retailer inspection policies on the distributor's preference for cold chain options. Transitioning from vaccine supply chain considerations, [27] highlighted the value of integrating Operational Research (OR) and AI models within Intelligent Decision Support Systems (IDSS). Abbass substantiated the collaborative integration approach using a meta-level modular representation utilizing constraint logic programming (CLP). The potential usefulness of integrated OR and AI is demonstrated through an illustration of a CLP-like meta-language applied to the dairy industry supply chain. Li et al. introduced an innovative combined ML and optimization-based framework to find a tradeoff between rescheduling frequency and the growing accumulation of delays. Also, they validated the importance of ML techniques and optimization algorithms for rescheduling decisions [28]. Using ML techniques like Artificial Neural Networks (ANN) and Non-dominated Sorting Genetic Algorithms (NSGA), Mouziraji et al.'s research aims to minimize geometrical deviations in wire cut machining. This study offered a novel perspective on studying geometric accuracy in wire cut drilling using a combination of ANN and NSGA, revealing optimal machining parameters [29].

Our paper presents contributions that bridge the gap between data science and optimization models through data-driven methodologies. In light of this, the quality of the product is expressed in terms of four important features, encompassing temperature, humidity, packaging types, and other environmental conditions, which are estimated by a synthetic dataset. After that, it plays a role of constraint in the model. The remainder of the paper is organized as follows. The problem definition, assumptions, and optimization are articulated in Section 3. The mathematical model, along with the framework of the dataset, is parented in the next section. Detailed information regarding the method of solving and the results of the study can be found in the following sections. Finally, in the last section, we offer some concluding remarks and suggest future research directions.

## 3 Problem Description

A warehouse generally stores a wide range of food items; and for each, a dedicated, environmentally controlled vault is operated. Without loss of generality, we address the cold storage management for an item. Our proposed model could easily be generalized to cover multiple items. We assume the storage of the warehouse and the number of orders are constrained, and it has employed a continuous review inventory control strategy. We further assume backordering is allowed and the lead-time demand for each item is a random variable. It is also worth mentioning that the other prominent aim is to determine different features of a product's quality, including temperature, humidity, packaging, and environmental conditions. To achieve this, a data-driven has been employed. In fact, we meticulously selected these four features for a regression model following a rigorous trial-and-error process. To develop a robust regression model for predicting product quality, a synthetic dataset comprising these essential features was generated using Python. The function of this regressing will be as a constraint in the presented model.

Additionally, to ensure a comprehensive analysis, the model is assembled based on several key assumptions:

- The quality of the product is assumed to fall within a predetermined value.
- Environmental conditions are categorized into three distinct modes, allowing for a more refined analysis of their impact on the cold chain.
- The demand quantity of the product is assumed to follow a uniform random distribution.

Considering these crucial aspects, our data-driven model aims to optimize warehouse management for the product, leading to decreasing costs and enhancing operational efficiency. The primary objective of this study is to minimize the total cost of our product, encompassing the total holding cost, fixed ordering cost, shortage cost, as well as costs associated with modifying temperature, humidity, packaging, and environmental conditions.

## 4 Model Formulation

In this section, we thoroughly present and elaborate on the proposed mathematical model, encompassing a detailed explanation of each model component. In the following, Table 1 displays the indices that are used in mathematical modeling. The model's parameters are shown in Table 2, and Table 3 showcases the decision variables within our model.

**Table 1** Index

| Index | Definition |
|---|---|
| $J$ | Index for 3 different packaging types and environmental conditions (j = 1,2,3) |

**Table 2** Parameters

| Parameter | Definition |
|---|---|
| $B$ | The average shortage of the product in year should be less than this amount |
| $Qual$ | The quality of product should be greater than this level |
| $f$ | The amount of warehouse space occupied by the product |
| $D$ | Demand quantity of the product |
| $A$ | Ordering cost of the product |
| $F$ | Available space for storing and holding the products at the warehouse |
| $h$ | Holding cost of one unit of the product |
| $\pi$ | Shortage cost per unit of the product |
| $e$ | Variable cost for temperature |
| $k$ | Variable cot for humidity |
| $b$ | Fixed cost of product's temperature distance to its ideal state |
| $d$ | Fixed cost of product's humidity distance to its ideal state |
| $\mu$ | Average demand for the product |
| $n$ | The number of product order shouldn't be greater than this amount |
| $M_j$ | The expenses associated with each packaging variant for the product |
| $N_j$ | The expenses associated with each environmental condition for the product |
| $Tu$ | Upper bound temperature for storing the product in warehouse satisfactorily |
| $Tl$ | Lower bound temperature for storing the product in warehouse properly |
| $HUu$ | Upper bound humidity for the storing the product in warehouse appropriately |
| $HUl$ | Lower bound humidity for the storing the product in warehouse suitably/nicely |
| $x1, x2, x3, x4$ | Linear Regression Function coefficients |
| $L$ | The constant number of Linear Regression Function |

**Table 3** Decision Variables

| Decision Variables | Definition |
|---|---|
| $Q$ | Ordering Amount (kg) of the product |
| $R$ | Reordering Point of the product (considered uniform distribution in the range of (100,200)) |
| $T$ | Temperature, the range considered between (-5,5) Celsius |
| $HU$ | Humidity, the range considered between (60,90) Percentage |
| $Y_j$ | 1, if packaging of type j is selected, 0 otherwise |
| $Z_j$ | 1, if environmental of type j is selected, 0 otherwise |

Based on the parameters and decision variables which have been defined in the aforementioned tables, the mathematical model will be:

$$Minimize\ TC = \left(\frac{D}{Q}\right)A + h\left(\frac{Q}{2} + R - \mu\right) + \left(\frac{D}{Q}\right) \times \pi \times \frac{(200-R)^2}{200} + \frac{D}{Q}(b - e(T - \text{Tu})) + \frac{D}{Q}(d - k(HU - \text{HUu})) + \sum_{j=1}^{3}\left(\frac{D}{Q}\right)M_j \times Y_j + \sum_{j=1}^{3}\left(\frac{D}{Q}\right)N_j \times Z_j \quad (1)$$

**S.t:**

$$x1(T) + x2(HU) + x3\left(\sum_{j=1}^{3} Y_j \times j\right) + x4\left(\sum_{j=1}^{3} Z_j \times j\right) + L \geq Qual \quad (2)$$

$$\left(\frac{D}{Q}\right) \times \left(\frac{(200-R)^2}{200}\right) \leq AB \quad (3)$$

$$f(Q + R) \leq F \quad (4)$$

$$\left(\frac{D}{Q}\right) \leq n \quad (5)$$

$$\sum_{j=1}^{3} Y_j = 1 \quad (6)$$

$$\sum_{j=1}^{3} Z_j = 1 \quad (7)$$

$$Hl \leq H \leq Hu \quad (8)$$

$$Tl \leq T \leq Tu \quad (9)$$

$$Rl \leq R \leq Ru \quad (10)$$

$$Q, R \geq 0 \quad (11)$$

$$Y_j, Z_j \in \{0,1\} \quad (12)$$

The Objective Function (OF) is declared in Eq. (1), which minimizes the total cost, including ordering, holding, shortage items, and adjusting the temperature, humidity, selecting modes of packaging, and environmental conditions. With regard to the shortage cost, the demand is assumed in the form of uniform distribution between 100 and 200. Eq. (2) represents the constraint which is related to product quality obtained by a data set. To create the dataset, synthetic data has been used to generate data, and data's features have four distinct variables: temperature, humidity, packaging mode, and environmental conditions. These attributes were thoughtfully chosen based on iterative experimentation, considering their significant influence on product quality. Data generation was achieved using the Python programming language with the help of the NumPy library. It takes the form of a Linear Regression Function, ensuring that the product quality exceeds the specified Qual amount. The Eq. (3) constraint implies that the average shortage of the product per year should be less than a certain limit (B). The constraint represented by Eq. (4) stipulates that the capacity occupied space of all products must not surpass the available storage capacity (F) designated for housing all products within the warehouse. Constraint (5) specifies that the frequency of orders for the desired product should remain below a predetermined threshold, denoted as 'n'. Eq. (6,7) serve to ensure that the model selects only one packaging type and environmental condition, optimally accommodating the product in the warehouse. The constraints of Eq. (8,9) establish boundaries for storage of the product in suitable temperature and humidity, ensuring the product remains within a safe range, free from decay or degradation. Lastly, the (11,12) constraints specify the sign and permissible values for decision variables in the model.

## 5 Solution Approach

The mathematical model is classified as a Mixed-Integer Non-Linear Programming (MINLP) model, which we solved using the Pyomo library in Python 3.8.8. The base amounts for each parameter in Table 4 are used to solve the model, and the output, which is the base result for this model, is declared in Table 5. As a matter of fact, the reason for using the Pyomo library was its efficiency and effectiveness for better use of sensitivity analysis in the later stages of the study. It is noteworthy to emphasize the results of this optimization are indeed optimal. All computations were executed on a Windows operating system, powered by an Intel(R) Core(TM) i7-10750H CPU @ 2.60GHz, 2592 Mhz, with 6 Cores and 12 Logical Processors.

**Table 4** Parameter amounts

| Parameter | Range |
|---|---|
| B | (0,6) |
| Qual | (40,100) |
| F | (1,3) |
| D | (5000,15000) |
| A | (100,500) |
| F | (1000,2000) |
| H | (12,24) |
| $\Pi$ | (4,10) |
| E | (10,20) |
| K | (10,20) |
| B | (100,400) |
| D | (100,400) |
| $\mu$ | 150 |
| N | (10,30) |
| Mj | {1:500, 2:1000, 3:1500} |
| Nj | {1:500, 2:1200, 3:2000} |
| Tu | 5 |
| Tl | -5 |
| Huu | 90 |
| HUl | 60 |
| x1, x2, x3, x4 | x1: -12.88, x2: -33.56, x3: 2.86, x4: 4.28 |
| L | 79.63 |

**Table 5** Output Values

| Objective Function | Decision Variables | | | |
|---|---|---|---|---|
| Total Cost | T | HU | Yj | Zj |
| 38444.8655 | 5 | 63.6 | Y1=1, Y2=Y3=0 | Z1=1, Z2=Z3=0 |

## 6  Discussion

Sensitivity analysis is a vital aspect of operations research, exploring how parameters affect the optimal solution. It delves into the optimal value of the objective function and the problem's best solution while taking into account significant parameter variations [30]. Sensitivity analysis aids in showcasing and forecasting the consequences of altering a model's input variables on the model's output. While certain parameters have a negligible impact on the model's outcomes, others exhibit a notable influence when subject to variations [31].In this approach, we evaluate the impact of altering a single parameter while assuming that the other parameters remain constant. If the optimal solution is found to be strongly influenced by adjustments in a parameter, we adjust its value and afterward rework the model. The essential goal is to determine a solution that indicates reduced sensitivity to parameter modifications. According to the parameters' value, sensitivity analysis is performed on parameters such as Qual, B, h, $\pi$, T, and HU.

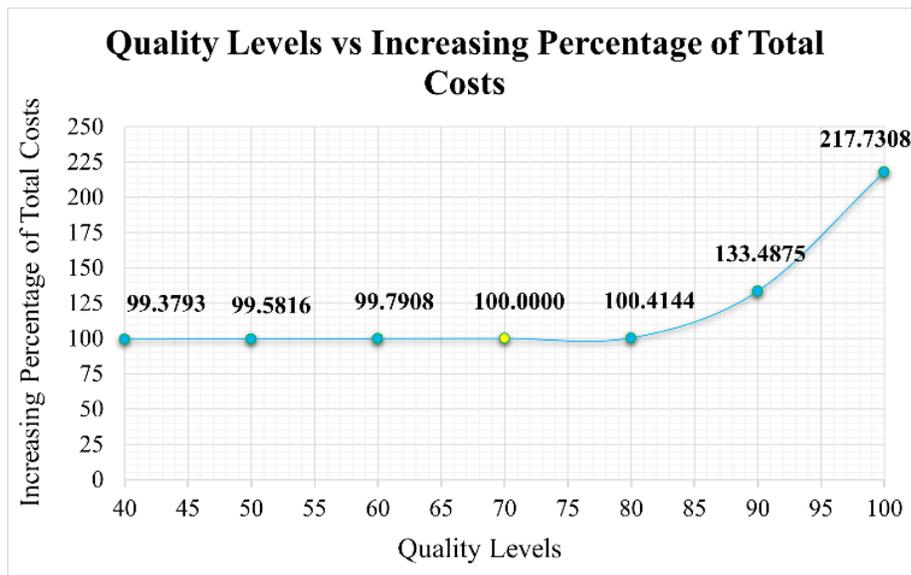

**Fig. 1** Quality levels vs. Total cost

The Fig. 1 illustrates the relation between quality and increasing the percentage of the total cost. Broadly speaking, the higher the quality is, the greater the increase in the percentage of the total cost will be, depicting an upward trend. A detailed analysis indicates that the increase in the percentage of total cost remains negligible from quality levels of 40% to 70%, and it results in the line chart seeming almost static in this range. Although there is a slight drop in the range of quality from 74% to 77%, the pattern subsequently regains its original trend. It's worth highlighting that the optimal quality level is at 80%, as this point allows us to achieve a reasonable quality without incurring a significant increase in costs. Managers would find it advantageous to target this quality level, given that the costs are relatively comparable to those at 40% quality. Notably, there is a marked surge in the quality from 80% to 100%. This growth is attributed to improved packaging and environmental conditions, both of which can influence both quality and total cost.

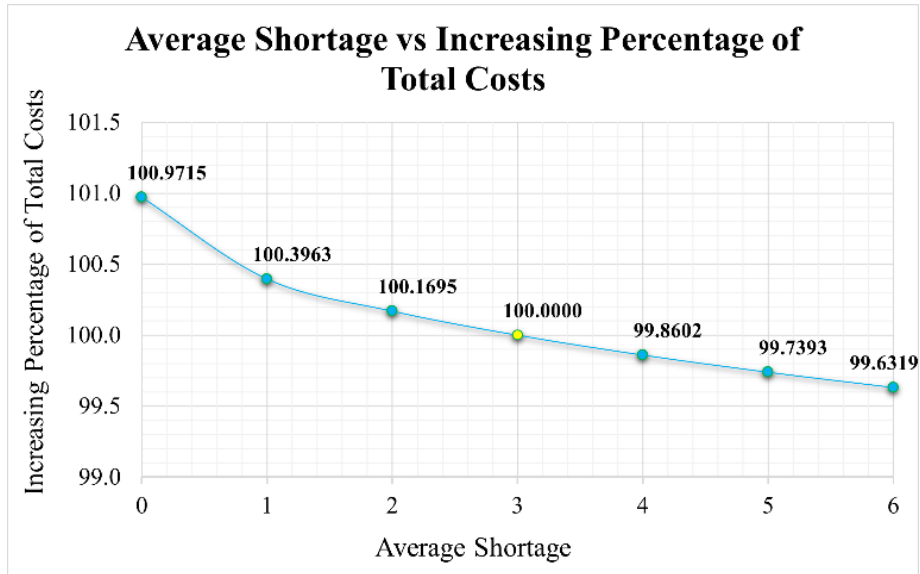

**Fig. 2** Average Shortage vs. Total cost

A notable correlation exists between the average shortage and the total cost, as can be seen in Fig. 2. As the average shortage cost rises, the total cost declines, reflecting an inverse relationship in a roughly rational function graph. However, as the average shortage cost increases, the fluctuation in total cost will reduce.

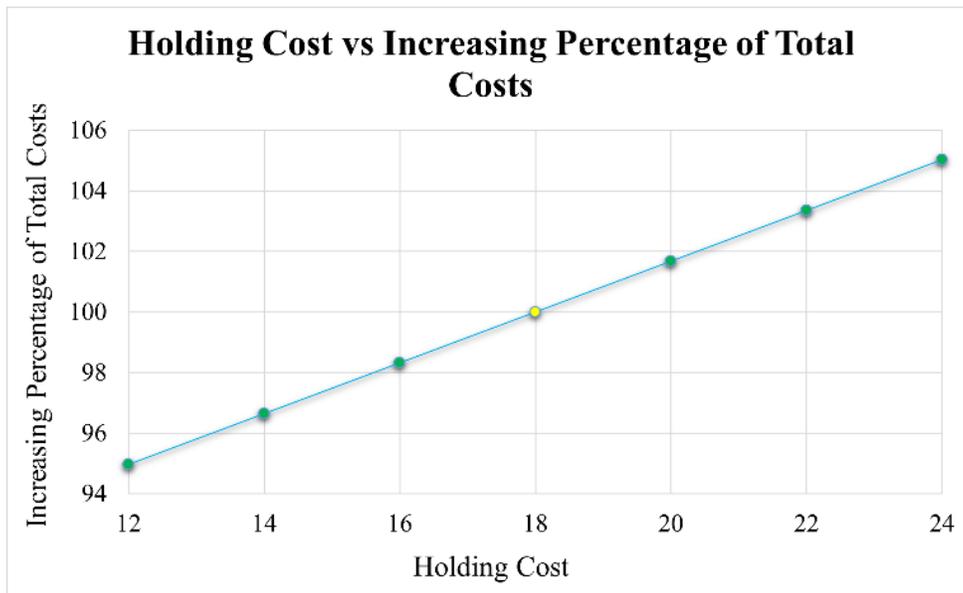

**Fig. 3** Holding Cost vs. Total cost

Holding costs constitute a significant component of the total costs, and any upward adjustment in holding cost with a fixed order size ordering system while keeping the other parameters constant results in a corresponding increase in the total cost. This relationship can be illustrated in Fig. 3.

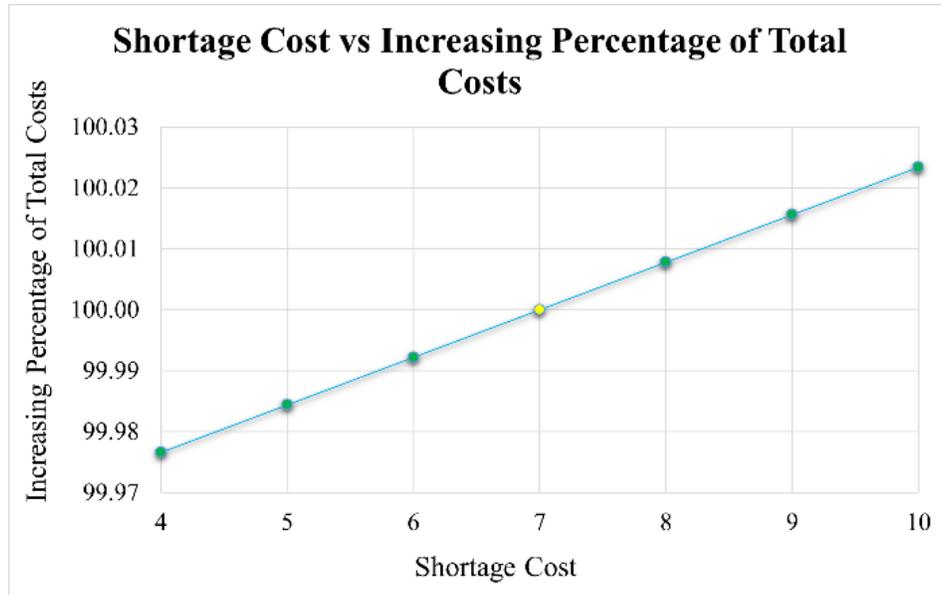

**Fig. 4** Shortage Cost vs. Total cost

In a fixed order size ordering system, increasing the shortage cost while keeping other factors constant can lead to an overall increase in the total cost. This relationship is visually represented in Fig. 4.

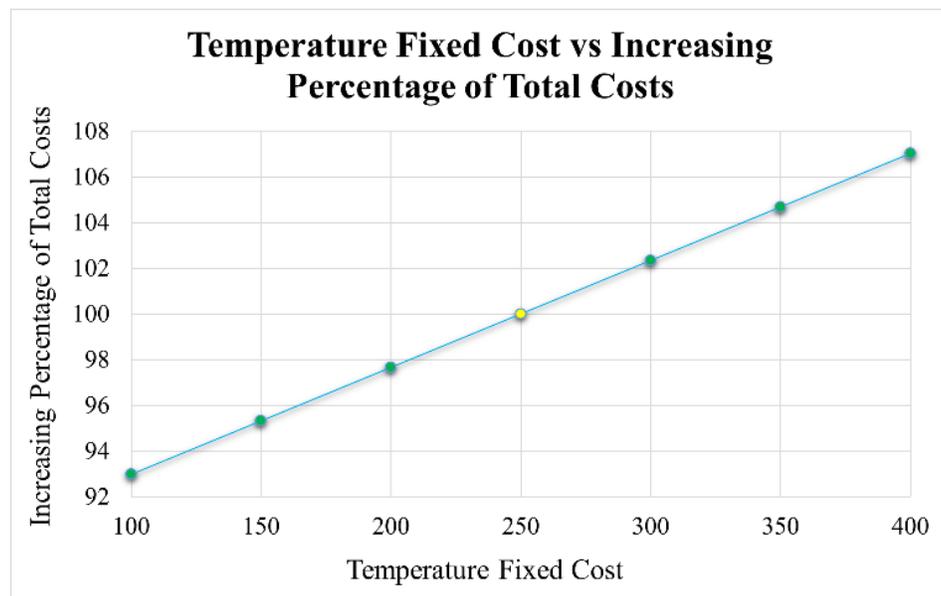

**Fig. 5** Temperature Cost vs. Total cost

The information extracted from Fig. 5 shows that in a fixed order size ordering system, an upsurge in the fixed costs associated with temperature while maintaining the other features static will lead to an upsurge in total costs.

## 7 Conclusion and Future Research

Our contribution to this study lies in the development of an optimization model coupled with an ML approach to ascertain the optimal warehouse features and find ordering system variables such as order quantities and reorder points. Moreover, these variables can be adjusted when altering the quality level to optimize customer demand. It is natural that achieving higher quality levels necessitates elevating the standards for temperature, humidity, and other environmental conditions.

By using the data being generated, a regression model has been constructed that establishes a relationship between product quality and four environmental conditions. This relationship was integrated into our non-linear mixed-integer programming model, enabling us to derive optimal solutions for these variables while minimizing costs. While the data being used was synthetic, it's important to note that this approach can readily transition to real-world scenarios, empowering decision-makers to enhance efficiency and make informed choices. Furthermore, our study employed sensitivity analysis to determine the optimal quality threshold with respect to total cost and explored the impact of varying parameters such as holding costs, shortage costs, humidity, and fixed temperature parameters on the total cost. The sensitivity analysis results are presented in Table 6.

**Table 6.** Summary Results

| Result of the sensitivity analysis on each parameter |
| --- |
| Increasing holding cost will result in higher total cost. |
| Raising the value of shortage cost leads to a proportional increase in the total cost. |
| Increasing average shortage leads to a decrease in total cost. |
| Higher quality results in higher costs. |
| Increasing humidity/temperature fix cost results in higher cost. |

Looking ahead, this study opens avenues for future exploration:

- Multi-Objective Model: This work can be extended to incorporate a multi-objective model that simultaneously considers product quality and costs, providing a more holistic approach to decision-making. In addition, other objectives encompassing environmental pollution and social considerations (enhancing the workforce) can be added to the model.

- Additional Quality Estimation Features: The model can be expanded to include additional features for estimating product quality, such as transportation mode, path distance, weather conditions, timestamps, and more, offering more comprehensive insights.

- Model Refinement: The model can be refined, and investigation between EOQ (Economic Order Quantity) and Periodic Review systems can provide a better intuition for decision-making and increase efficiency. Additionally, the inflation rate can be considered in the model to make it more realistic.

- Solving Method: Different solving methods, such as Meta-Heuristic, can be defined, especially for large-scale datasets. For example, the model can be solved using Genetic algorithms and particle swarm optimization.

- Parameters: The model can be equipped with holding and ordering costs with discounts, enhancing real-world applicability. Additionally, defining parameters in a random way can better capture the unpredictability of the real world.